\documentclass{amsart}
\usepackage{graphicx}
\newtheorem{thm}{Theorem}[section]

\newtheorem{lem}[thm]{Lemma}

\newtheorem{property}[thm]{Property}
\theoremstyle{definition}
\newtheorem{defn}[thm]{Definition}
\theoremstyle{remark}

\numberwithin{equation}{section}

\begin{document}

\title[LOSSES IN $M/GI/m/n$ QUEUES]{LOSSES IN $M/GI/m/n$ QUEUES}%
\author{Vyacheslav M. Abramov}%
\address{School of Mathematical Sciences, Monash University, Building 28M,
Clayton Campus, Wellington road, VIC-3800, Australia}%
\email{vyacheslav.abramov@sci.monash.edu.au}%


\begin{abstract}
The $M/GI/m/n$ queueing system with $m$ homogeneous servers and
the finite number $n$ of waiting spaces is studied. Let $\lambda$
be the customers arrival rate, and let $\mu$ be the reciprocal of
the expected service time of a customer. Under the assumption $\lambda=m\mu$
it is proved that the expected number of
losses during a busy period is the same  value for all $n\geq1$,
while in the particular case of the Markovian system $M/M/m/n$ the
expected number of losses during a busy period is $\frac{m^m}{m!}$
for all $n\geq0$. Under the additional assumption that the
probability distribution function of a service time belongs to the
class NBU or NWU, the paper establishes simple inequalities for
those expected numbers of losses in $M/GI/m/n$ queueing systems.
\end{abstract}
\maketitle
\section{Introduction}
\noindent Analysis of loss queueing systems is very important from
both the theoretical and practical points of view.
While the multiserver loss queueing system $M/GI/m/0$ and its network
extensions have been intensively studied (see the review paper of Kelly
\cite{Kelly - review}, the book of Ross \cite{Ross - book} and
references in these sources),
 the information about $M/GI/m/n$ queueing
systems ($n\geq 1$) is very scanty, because explicit results for
characteristics of these queueing systems are unknown. (In
the present paper, for multiserver queueing systems the notation
$M/GI/m/n$ is used, where $m$ denotes the number of servers and $n$
denotes the number of waiting places. Another notation which is also
acceptable in the literature is $M/GI/m/m+n$.)

From the practical point of view, $M/GI/m/n$ queueing systems serve
as a model for telephone systems, where $n$ is the maximally
possible number of calls that can wait in the line before their
service start. The loss probability is one of the most significant
performance characteristics. In the present paper, we study the
expected number of losses during a busy period (the characteristic
closely related to the stationary loss probability) under the
assumption that the arrival rate ($\lambda$) is equal to the maximum
service capacity ($m\mu$), which seems to be the most interesting
from the theoretical point of view.

There are two main reasons for studying this case.

The first reason is that the case $\lambda=m\mu$ is a critical case for
queueing systems with $m$ identical servers, i.e. the case associated
with critically loaded systems. The theoretical and practical interest
in studying heavily loaded loss systems is very high, and there are many
results in the literature related to the analysis of the loss probability
in heavily loaded systems. The asymptotic results for losses in heavily
loaded single server systems ($n\to\infty$) such as $M/GI/1/n$ and $GI/M/1/n$
and for associated models of telecommunication systems and dams have been studied
in \cite{Abramov 1997}, \cite{Abramov 2002}, \cite{Abramov 2004}, \cite{Abramov 2007},
\cite{Abramov 2008b} and \cite{Whitt 2004a}.
Heavy-traffic analysis of losses in heavily loaded multiserver systems have been provided
in \cite{Abramov 2007b}, \cite{Whitt 2004a}, \cite{Whitt 2004b} and \cite{Whitt 2005}.
The mathematical foundation of heavy traffic theory can be found in the
textbook of Whitt \cite{Whitt-book}.
Although the case $\lambda=m\mu$ is idealistic, it enables us to understand the possible
behaviour of the system in certain cases when the values $\lambda$ and $m\mu$ are close
and approach one another as $n$ increases to infinity. (Obtaining nontrivial results in the cases
$\lambda<m\mu$ and $\lambda>m\mu$ is a hard problem, so the analytic investigation of the
aforementioned asymptotic behaviour as $n$ increases to infinity is difficult.)

The second reason is that $\lambda=m\mu$ is an interesting
theoretical case associated with an extension of the following
non-trivial property of the symmetric random walk. Let $X_1$, $X_2$,
\ldots, $X_i$, \ldots, be a sequence of independent and identically
distributed random variables taking the values $\pm1$ with the equal
probability $\frac{1}{2}$. Let $S_0=0$, and $S_{i+1}=S_i+X_{i+1}$,
$i\geq0$, be a symmetric random walk, and let $t=\tau$ be the first
time instant after $t=0$ when this random walk returns to zero, i.e.
$S_\tau=0$. It is known that the expected number of level-crossings
through any level $n\geq1$ (or $n\leq-1$) is equal to $\frac{1}{2}$
independently of that level. The mentioning of this fact (but in a
slightly different formulation) can be found in Szekely
\cite{Szekely 1986}, and its proof is given in Wolff \cite{Wolff
1989}, p.411. The reformulation of this fact in terms of queueing
theory is as follows. Consider $M/M/1/n$ queueing system with equal
arrival and service rates. For this system, the expected number of
losses during a busy period is equal to 1 for all $n\geq0$. It has
been recently noticed that this property holds true for $M/GI/1/n$
queueing systems. Namely, it was shown in several recent papers (see
Abramov \cite{Abramov 1991a}, \cite{Abramov 1991b}, \cite{Abramov
1997}, Righter \cite{Righter 1999}, Wolff \cite{Wolff 2002}), that
under mutually equal expectations of interarrival and service time,
the expected number of losses during a busy period is equal to 1 for
all $n\ge 0$. Further extension of this property to queueing systems
with batch arrivals have been given in Abramov \cite{Abramov 2001a},
Wolff \cite{Wolff 2002} and Pek\"oz, Righter and Xia  \cite{Pekoz et
al 2003}. Applications of the aforementioned property of losses can
be found in \cite{Abramov 2004} for analysis of lost messages in
telecommunication systems and in \cite{Abramov 2007} for optimal
control of large dams. Further relevant results associated with the
properties of losses have been obtained in the paper by  Pek\"oz,
Righter and Xia \cite{Pekoz et al 2003}. They solved a
characterization problem associated with the properties of losses in
$GI/M/1/n$ queues and established similar properties for $M/M/m/n$
and $M^X/M/m/n$ queueing systems. Recently, a similar property
related to consecutive losses in busy periods of $M/GI/1/n$ queueing
systems has been discussed in \cite{Abramov 2009a}. It follows from
the results obtained in this paper that for $M/GI/1/n$ queueing
systems with mutually equal expectations of interarrival and service
times, the expected number of losses in series containing at least
$k>1$ consecutive losses during a busy period generally depends on
$n$. However, for $M/M/1/n$ queueing systems with equal arrival and
service rates that expected number of consecutive losses during a
busy period is the same constant (depending on the value $k$) for
all $n\geq0$.

The aim of the present paper is further theoretical contribution to
this theory of losses, now to the theory of multiserver loss
queueing systems. On the basis of the aforementioned results on
losses in $M/GI/1/n$ and $M/M/m/n$ queueing systems we address the
following open question. \textit{Does the result on losses in}
$M/M/m/n$ \textit{queueing systems remain true for those} $M/GI/m/n$
\textit{too?}

The answer on this question is not elementary.
On one hand,
under the assumption $\lambda=m\mu$ the expected numbers of losses in $M/GI/m/0$ and $M/GI/m/n$ queueing systems ($m\geq2$ and $n\geq1$) during their busy periods are different. A simple example for this confirmation can be built
for $M/GI/2/1$ queueing systems having the service time
distribution $G(x)=1-p\mbox{e}^{-\mu_1x}-q\mbox{e}^{-\mu_2x}$,
$p+q=1$. The analysis of the stationary characteristics for these
systems, resulting in an analysis of losses during a busy period,
can be provided explicitly. Specifically, the structure of the
$9\times9$ Markov chain intensity matrix for the states of the
Markov chain associated with an $M/GI/2/1$ queueing system shows a
clear difference between the structure of the stationary
probabilities in $M/GI/2/1$ queues and that in $M/GI/2/0$ queues
given by the Erlang-Sevastyanov formulae. So, the parameters $p$,
$q$, $\mu_1$ and $\mu_2$ can be chosen such that the expected
number of losses during busy periods in these two queueing systems
will be different.

On the other hand, the property of losses, which is similar to the aforementioned one, indeed holds. The correctness of this similar property for multiserver $M/GI/m/n$ queueing systems is proved in the present paper. Namely, we establish the following results.

Let $L_{m,n}$ denote the number of losses
during a busy period of the $M/GI/m/n$ queueing system, let
$\lambda$, $\mu$ be the arrival rate and, respectively, the
reciprocal of the expected service time, and let $m$, $n$ denote
the number of servers and, respectively, the number of waiting
places. We will prove that, under the assumption $\lambda=m\mu$,
the expected number of losses during a busy period of the $M/GI/m/n$
queueing system, $\mathrm{E}L_{m,n}$, is the same for all
$n\geq1$, which is \textit{not} generally the same as that for the
$M/GI/m/0$ loss queueing system (when $n=0$). In addition, if the probability
distribution function of the service time belongs to the class NBU
(New Better than Used), then $\mathrm{E}L_{m,n}= \frac{cm^m}{m!}$,
where a constant $c\geq1$ is independent of $n\geq1$. In the
opposite case of the NWU (New Worse than Used) service time
distribution we correspondingly have
$\mathrm{E}L_{m,n}=\frac{cm^m}{m!}$ with a constant $c\leq1$
independent of $n\geq1$ as well. (The constant $c$ becomes equal
to 1 in the case of exponentially distributed service times.)
Recall that a probability distribution function $\Xi(x)$ of a
nonnegative random variable is said to belong to the class NBU if
for all $x\geq0$ and $y\geq0$ we have
$\overline{\Xi}(x+y)\leq\overline{\Xi}(x)\overline{\Xi}(y)$, where
$\overline{\Xi}(x)=1-{\Xi}(x)$. If the opposite inequality holds,
i.e. $\overline{\Xi}(x+y)\geq\overline{\Xi}(x)\overline{\Xi}(y)$,
then $\Xi(x)$ is said to belong to the class NWU.

The proof of the main results of this paper is based on an
application of the level-crossing approach to the special type
stationary processes. The construction of the level-crossings
approach used in this paper is a substantially extended version of that used in the earlier
papers by the author (e.g. \cite{Abramov 1991a}, \cite{Abramov
1994}, \cite{Abramov 2001b}, \cite{Abramov 2001c}, \cite{Abramov
2006} and \cite{Abramov 2008}) and by Pechinkin
\cite{Pechinkin 1987}. It uses modern geometric methods of analysis and involves an algebraically close system of processes and a nontrivial construction of deleting intervals and merging the ends together with nontrivial applications of the PASTA property.

Throughout the paper, it is assumed that $m\geq 2$. (This is not
the loss of generality since the case $m=1$ is known, see
\cite{Abramov 1997}, \cite{Righter 1999} and \cite{Wolff 2002}.)

The paper is organized as follows. In Section \ref{Markovian}, which
is the first part of the paper, $M/M/m/n$ queueing systems are
studied. The results for $M/M/m/n$ queueing systems are then used in
Section \ref{Non-Markovian}, which is the second part of the paper,
in order to study $M/GI/m/n$ queueing systems. The study in both of
Sections \ref{Markovian} and \ref{Non-Markovian} is based on the
level-crossing approach. The construction of level-crossings for
$M/M/m/n$ queueing systems is then developed for $M/GI/m/n$ queueing
systems as follows. The stationary processes associated with these
queueing systems is considered, and the stochastic relations between
the times spent in state $m-1$ associated with $m-1$ busy servers
during a busy period of $M/GI/m/n$ ($n\geq1$) and $M/GI/m-1/0$
queueing systems are established. To prove these stochastic
relations, some ideas from the paper of Pechinkin \cite{Pechinkin
1987} are involved to adapt and develop the level-crossing method
for the problems of the present paper. The obtained stochastic
relations are crucial, and they are then used to prove the main
results of the paper in Section 4. In Section 5, possible
development of the results for $M^X/GI/m/n$ queueing systems with
batch arrivals is discussed.

\section{The $M/M/m/n$ queueing system}\label{Markovian}
\noindent In this section, the Markovian $M/M/m/n$ loss
queueing system is studied with the aid of the level-crossings approach, in
order to establish some relevant properties of this queueing
system. Those properties are then developed for $M/GI/m/n$
queueing systems in the following sections.

Let $f(j)$, $1\leq j\leq n+m+1$, denote the number of customers
arriving during a busy period who, upon their arrival, meet $j-1$
customers in the system. It is clear that $f(1)=1$ with
probability 1. Let $t_{j,1}$, $t_{j,2}$,\ldots, $t_{j,f(j)}$ be
the instants of arrival of these $f(j)$ customers, and let
$s_{j,1}$, $s_{j,2}$,\ldots, $s_{j,f(j)}$ be the instants of the
service completions when there remain only $j-1$ customers in the
system. Notice, that $t_{n+m+1,k}=s_{n+m+1,k}$ for all
$k=1,2,\ldots,f(n+m+1)$.

For $1\leq j\leq n+m$ let us consider the intervals
\begin{equation}\label{2.1}
(t_{j,1}, s_{j,1}], (t_{j,2}, s_{j,2}],\dots, (t_{j,f(j)},
s_{j,f(j)}].
\end{equation}
Then, by incrementing index $j$ we have the following intervals
\begin{equation}\label{2.2}
(t_{j+1,1}, s_{j+1,1}], (t_{j+1,2}, s_{j+1,2}],\ldots,
(t_{j+1,f(j+1)}, s_{j+1,f(j+1)}].
\end{equation}
Delete the intervals of \eqref{2.2} from those of \eqref{2.1} and
merge the ends, that is each point $t_{j+1,k}$ with the
corresponding point $s_{j+1,k}$, $k=1,2,...,f(j+1)$ (see Figure 1).

\begin{figure}
\includegraphics[width=15cm, height=20cm]{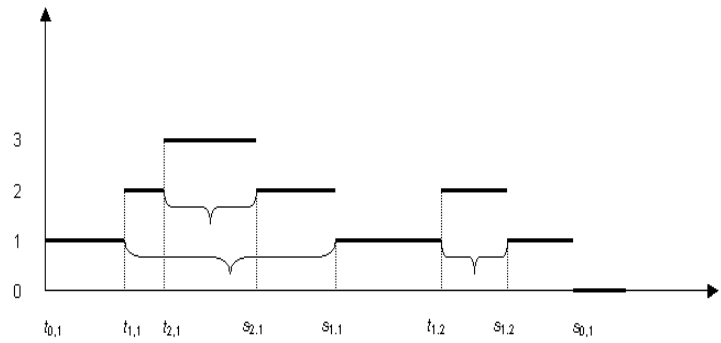}
\caption{Level crossings during a busy period in a Markovian system}
\end{figure}

Then $f(j+1)$ has the following properties. According to the
property of the lack of memory of the exponential distribution, the
residual service time for a service completion, after the procedure of deleting the interval
and merging the ends as it is indicated above, remains exponentially
distributed with parameter $\mu \min(j,m)$. Therefore, the number of
points generated by merging the ends within the given interval
$(t_{j,1}, s_{j,1}]$ coincides in distribution with the number of
arrivals of the Poisson process with rate $\lambda$ during an
exponentially distributed service time with parameter $\mu
\min(j,m)$. Namely, for $1\leq j\leq m-1$ we obtain
$$
\mathrm{E}\{f(j+1)|f(j)=1\}=\sum_{u=1}^\infty
u\int_0^\infty\mbox{e}^{-\lambda x}\frac{(\lambda x)^u}{u!}j\mu
\mbox{e}^{-j\mu x}\mbox{d}x=\frac{\lambda}{j\mu}.
$$
Considering now a random number $f(j)$ of intervals (\ref{2.1}) we
have
\begin{equation}\label{2.4}
\mathrm{E}\{f(j+1)|f(j)\}=\frac{\lambda}{j\mu}f(j).
\end{equation}
Analogously, denoting the load of the system by
$\rho=\frac{\lambda}{m\mu}$, for $m\leq j\leq m+n$ we have
\begin{equation}\label{2.5}
\mathrm{E}\{f(j+1)|f(j)\}=\frac{\lambda}{m\mu}f(j)=\rho f(j).
\end{equation}
The properties \eqref{2.4} and \eqref{2.5} mean that the
stochastic sequence
\begin{equation}\label{2.5+}
\left\{f(j+1)\left(\frac{\mu}{\lambda}\right)^j\prod_{i=1}^{j}\min(i,m),\mathcal{F}_{j+1}\right\},\
\ \mathcal{F}_j=\sigma\{f(1), f(2), \ldots, f(j)\},
\end{equation}
forms a martingale.

It follows from \eqref{2.5+} that for $0\leq j\leq m-1$
\begin{equation}\label{2.6}
\mathrm{E}f(j+1)=\frac{\lambda^j}{j!\mu^j},
\end{equation}
and for $m\leq j\leq m+n$
\begin{equation}\label{2.7}
\mathrm{E}f(j+1)=\frac{\lambda^m}{m!\mu^m}~\rho^{j-m}.
\end{equation}
For example, when $\rho=1$ from (\ref{2.7}) we obtain the
particular case of the result of Pek\"oz, Righter and Xia
\cite{Pekoz et al 2003}:
$\mathrm{E}L_{m,n}$=$\mathrm{E}f(n+m+1)=\frac{m^m}{m!}$ for all
$n\ge 0$, where $L_{m,n}$ denotes the number of losses during a
busy period of the $M/M/m/n$ queueing system.

\smallskip

Next, let $B(j)$ be the period of time during a busy cycle of the
$M/M/m/n$ queueing system when there are exactly $j$ customers in
the system. For $0\leq j\leq m+n$ we have:
\begin{equation}\label{2.8}
\lambda\mathrm{E}B(j)=\mathrm{E}f(j+1)=\begin{cases}\frac{\lambda^{j}}{j!\mu^{j}}, &\mbox{for} \ 0\leq j\leq m-1,\\
\frac{\lambda^{m}}{m!\mu^{m}}\rho^{j-m},&\mbox{for} \ m\leq j\leq m+n.
\end{cases}
\end{equation}
Now, introduce the following notation. Let $T_{m,n}$ denote the length of a busy
period of the $M/M/m/n$ queueing system, let $T_{m,0}$ denote the length of a
busy period of the $M/M/m/0$ queueing system with the same arrival
and service rates as in the initial $M/M/m/n$ queueing system, and
let $\zeta_n$ denote the length of a busy period of $M/M/1/n$ queueing system
with arrival rate $\lambda$ and service rate $\mu m$. From
(\ref{2.6})-(\ref{2.8}) for the expectation of a busy period of
the $M/M/m/n$ queueing system we have
\begin{equation}\label{2.10}
\mathrm{E}T_{m,n}=\sum_{j=1}^{n+m}\mathrm{E}B(j)=\sum_{j=1}^{m-1}
\frac{\lambda^{j-1}}{j!\mu^{j}}+\frac{\lambda^{m-1}}{m!\mu^{m}}\sum_{j=0}^{n}\rho^j.
\end{equation}
In turn, for the expectation of a busy period of the $M/M/m/0$
queueing system we have
\begin{equation}\label{2.11}
\mathrm{E}T_{m,0}=\sum_{j=1}^{m}\mathrm{E}B(j)= \sum_{j=1}^{m}
\frac{\lambda^{j-1}}{j!\mu^{j}},
\end{equation}
where (\ref{2.11}) is the particular case of (\ref{2.10}) where
$n=0$.

It is clear that $T_{m,n}$ contains one busy period
$T_{m-1,0}$, where the subscript $m-1$ underlines that there are
$m-1$ servers, and a random number of independent busy periods,
which will be called \textit{orbital busy periods}. Denote an
orbital busy period by $\zeta_n$. (It is assumed that an orbital
busy period $\zeta_n$ starts at instant when an arriving customer
finds $m-1$ servers busy and occupies the $m$th server, and it
finishes at the instant when after a service completion there at
the first time remain only $m-1$ busy servers.) Therefore,
denoting the independent sequence of identically distributed
orbital busy periods by $\zeta_n^{(1)}$, $\zeta_n^{(2)}$,..., we
have
\begin{equation}\label{2.12}
T_{m,n}{\buildrel{d}\over =}T_{m-1,0}+\sum_{i=1}^\kappa\zeta_n^{(i)},
\end{equation}
where $\kappa$ is the random number of the aforementioned orbital
busy periods and ${\buildrel d\over =}$ means an equality in distribution. It follows from (\ref{2.10}), (\ref{2.11}) and
(\ref{2.12})
\begin{equation}\label{2.13}
\mathrm{E}\sum_{i=1}^\kappa\zeta_n^{(i)}=\frac{\lambda^{m-1}}
{m!\mu^{m}}\sum_{j=0}^{n}\rho^j.
\end{equation}
On the other hand, the expectation of an orbital busy period
$\zeta_n$ is
$$\mathrm{E}\zeta_n=\frac{1}{m\mu
}\sum_{j=0}^{n}\rho^j$$
(this can be easily checked, for example,
by the level-crossings method \cite{Abramov 1991a}, \cite{Abramov
2001b} and an application of Wald's identity \cite{Feller 1966},
p. 384), and we obtain
\begin{equation}\label{2.14}
\mathrm{E}\kappa=\frac{\lambda^{m-1}}{(m-1)!\mu^{m-1}}.
\end{equation}
Thus, $\mathrm{E}\kappa$ coincides with the expectation of the
number of losses during a busy period in the $M/M/m/0$ queueing
system. In the case $\rho=1$ we have
$\mathrm{E}\kappa=\frac{m^m}{m!}$.

\section{$M/GI/m/n$ queueing systems}\label{Non-Markovian}
In this section, the inequalities between the times spent in the
state $m-1$ in the $M/GI/m/n$ ($n\geq1$) and $M/GI/m/0$ queueing
systems during their busy periods are derived.

 Consider two queueing systems: $M/GI/m/n$ ($n\geq1$)
and $M/GI/m/0$ both having the same arrival rate $\lambda$ and
probability distribution function of a service time $G(x)$,
$\frac{1}{\mu}=\int_0^\infty x\mbox{d}G(x)<\infty$. Let
$T_{m,n}(m-1)$ denote the time spent in the state $m-1$  during
its busy period (i.e. the total time during a busy period when
$m-1$ servers are occupied) of the $M/GI/m/n$ queueing system, and
let $T_{m,0}(m-1)$ have the same meaning for the $M/GI/m/0$
queueing system.

We prove the following lemma.

\begin{lem}\label{Lemma3.1} Under the assumption that the service time
distribution $G(x)$ belongs to the class NBU (NWU),
\begin{equation}\label{3.1}
T_{m,n}(m-1)\geq_{st} (\mbox{resp.} \ \leq_{st})\ T_{m,0}(m-1),
\end{equation}
\end{lem}

\begin{proof}. The proof of the lemma is relatively long. In order to make
it transparent and easily readable we strongly indicate the steps of
this proof given by several propositions (properties). There are
also six figures (Figures 2-7) illustrating the constructions in the
proof. Each of these figures contain two graphs. The first (upper)
of them indicates the initial (or intermediate) possible path of the
process (sometimes two-dimensional), while the second (lower) one
indicates the part of the path of one or two-dimensional process
after a time scaling or specific transformation (e.g. in Figure 5).
Arc braces in the graphs indicate the intervals that should be
deleted and their ends merged.

Two-dimensional processes are shown as parallel graphs. For example,
there are two parallel processes in Figure 3 which are shown in the
upper graph, and there are two parallel processes which are shown in
the lower graph. The same is in Figures 4, 6 and 7.

\smallskip
For the purpose of the present paper we use strictly stationary
processes of order 1 or \textit{strictly 1-stationary processes}.
Recall the definition of a strictly stationary process of order
$n$ (see \cite{PlanetMath}, p.206).

\begin{defn}\label{defn1}
The process $\xi(t)$ is said to be \textit{strictly stationary of
order $n$} or \textit{strictly $n$-stationary}, if for a given
positive $n<\infty$, any $h$ and $t_1$, $t_2$,\ldots, $t_n$ the
random vectors
\begin{equation*}
\Big(\xi(t_1), \xi(t_2),\ldots, \xi(t_n)\Big) \ \mbox{and} \
\Big(\xi(t_1+h), \xi(t_2+h),\ldots, \xi(t_n+h)\Big)
\end{equation*}
have identical joint distributions.
\end{defn}

If $n=1$ then we have strictly 1-stationary processes satisfying
the property:
\begin{equation*}
\mathrm{P}\{\xi(t)\leq x\}=\mathrm{P}\{\xi(t+h)\leq x\}.
\end{equation*}
The probability distribution function $\mathrm{P}\{\xi(t)\leq x\}$
in this case will be called \textit{limiting stationary
distribution}.

The class of strictly 1-stationary processes is wider than the
class of strictly stationary processes, where it is required that
for all finite dimensional distributions
\begin{equation*}
\mathrm{P}\{(\xi(t_1), \xi(t_2),\ldots,\xi(t_k))\in
B_k\}=\mathrm{P}\{(\xi(t_1+h), \xi(t_2+h),\ldots,\xi(t_k+h))\in
B_k\},
\end{equation*}
for any $h$ and any Borel set $B_k\subset\mathbb{R}^k$. The reason
of using strictly 1-stationary processes rather than strictly
stationary processes themselves is that, the operation of deleting
intervals and merging the ends is algebraically close with respect
to strictly 1-stationary processes, and it is \textit{not} closed
with respect to strictly stationary processes. The last means that
if $\xi(t)$ is a strictly 1-stationary process, then for any $h>0$
and arbitrary $t_0$ the new process
$$
\xi_1(t)=\begin{cases}\xi(t), &\mbox{if} \ t\leq t_0,\\
\xi(t+h), &\mbox{if} \ t>t_0
\end{cases}
$$
is also strictly 1-stationary and has the same one-dimensional
distribution as the original process $\xi(t)$. The similar
property is not longer valid for strictly stationary processes. If
$\xi(t)$ is a strictly stationary process, then generally
$\xi_1(t)$ is not strictly stationary.

In the following the prefix `strictly' will be omitted, so
strictly stationary and strictly 1-stationary processes will be
correspondingly called stationary and 1-stationary processes.

\smallskip
Let us introduce $m$ independent and identically distributed stationary renewal
processes (denoted below ${\bf x}_m(t)$) with a renewal period
having the probability distribution function $G(x)$.

On the basis of these renewal processes we build the stationary
$m$-dimensional Markov process ${\bf x}_m(t)=\{\xi_1(t),
\xi_2(t),\ldots,\xi_m(t)\}$, the coordinates $\xi_k(t)$,
$k=$1,2,\ldots, $m$ of which are the residual times to the next
renewal times in time moment $t$, following in an ascending order.

\smallskip
Let us now consider the two $(m+1)$-dimensional Markov
processes corresponding to the $M/GI/m/n$ ($n\geq1$) and
$M/GI/m/0$ queueing systems, which are denoted by ${\bf
y}_{m,n}(t)$ and ${\bf y}_{m,0}(t)$.
Let $Q_{m,n}(t)$ denote the stationary queue-length process  (the
number of customers in the system) of the $M/GI/m/n$ queueing system, and
let $Q_{m,0}(t)$ denote the stationary queue-length process
corresponding to the $M/GI/m/0$ queueing system. We have:
\begin{equation*}\label{3.3}
{\bf y}_{m,n}(t)=\left\{\eta_1^{(m,n)}(t),
\eta_2^{(m,n)}(t),\ldots, \eta_{m}^{(m,n)}(t), Q_{m,n}(t)\right\},
\end{equation*}
where
$$
\left\{\eta_{m-\nu_{m,n}(t)+1}^{(m,n)}(t),
\eta_{m-\nu_{m,n}(t)+2}^{(m,n)}(t),\ldots,
\eta_{m}^{(m,n)}(t)\right\}
$$
are the ordered residual service times corresponding to
$\nu_{m,n}(t)=\min\{m, Q_{m,n}(t)\}$ customers in service in time
$t$, and
$$
\left\{\eta_{1}^{(m,n)}(t), \eta_{2}^{(m,n)}(t),\ldots,
\eta_{m-\nu_{m,n}(t)}^{(m,n)}(t)\right\}=\{{0,0,\ldots,0}\}
$$
all are zeros. Analogously,
\begin{equation*}\label{3.6}
{\bf y}_{m,0}(t)=\left\{\eta_1^{(m,0)}(t),
\eta_2^{(m,0)}(t),\ldots, \eta_{m}^{(m,0)}(t), Q_{m,0}(t)\right\},
\end{equation*}
only replacing the index $n$ with 0.

Let us delete all time intervals of the process ${\bf y}_{m,n}(t)$
related to the $M/GI/m/n$ queueing system ($n\geq1$) where there
are more than $m-1$ or less than $m-1$ customers and merge the
ends. Remove the last component of the obtained process which is
trivially equal to $m-1$. Then we get the new ($m-1$)-dimensional
Markov process (in the following the prefix `Markov' will be
omitted and only used in the places where it is meaningful):
\begin{equation*}\label{3.7}
\widehat{{\bf y}}_{m-1,n}(t)=\left\{\widehat{\eta}_1^{(m,n)}(t),
\widehat{\eta}_2^{(m,n)}(t),\ldots,
\widehat{\eta}_{m-1}^{(m,n)}(t)\right\},
\end{equation*}
the components of which are now denoted by hat. All of the
components of this vector are 1-stationary, which is a consequence
of the existence of the limiting stationary probabilities of the
processes $\eta_j^{(m,n)}(t)$, $j=1,2,\ldots,m$ (e.g. Tak\'acs
\cite{Takacs 1969}) and consequently those of the processes
$\widehat{\eta}_j^{(m,n)}(t)$, $j=1,2,\ldots,m$. The joint
limiting stationary distribution of the process $\widehat{{\bf
y}}_{m-1,n}(t)$ can be obtained by conditioning of that of the
processes ${{\bf y}}_{m,n}(t)$ given $Q_{m,n}(t)=m-1$.

The similar operation of deleting intervals and merging the ends,
where there are less than $m-1$ customers in the system, for the
process ${\bf y}_{m-1,0}(t)$ is used. We correspondingly have
\begin{equation*}\label{3.8}
\widehat{{\bf y}}_{m-1,0}(t)=\left\{\widehat{\eta}_1^{(m,0)}(t),
\widehat{\eta}_2^{(m,0)}(t),\ldots,
\widehat{\eta}_{m-1}^{(m,0)}(t)\right\}.
\end{equation*}

\smallskip
We establish the following elementary property related to the
$M/GI/1/n$ queueing systems, $n$=0,1,\ldots
\begin{property}\label{prop1}
\begin{equation}\label{fact-1}
\mathrm{P}\left\{\eta_1^{(1,n)}(t)\in
B_1~|~Q_{1,n}(t)\geq1\right\}=\mathrm{P}\{{\bf x}_1(t)\in B_1\},
\end{equation}
for any Borel set $B_1\subset\mathbb{R}^1$.
\end{property}

\smallskip
\begin{proof}
Delete all of the intervals where the server is free and merge the
corresponding ends (see Figure 2). Then in the new time scale, the
processes all are structured as a stationary renewal process with
the length of a period having the probability distribution function
$G(x)$. Therefore \eqref{fact-1} follows.
\end{proof}

\begin{figure}\label{fig1}
\includegraphics[width=15cm,height=20cm]{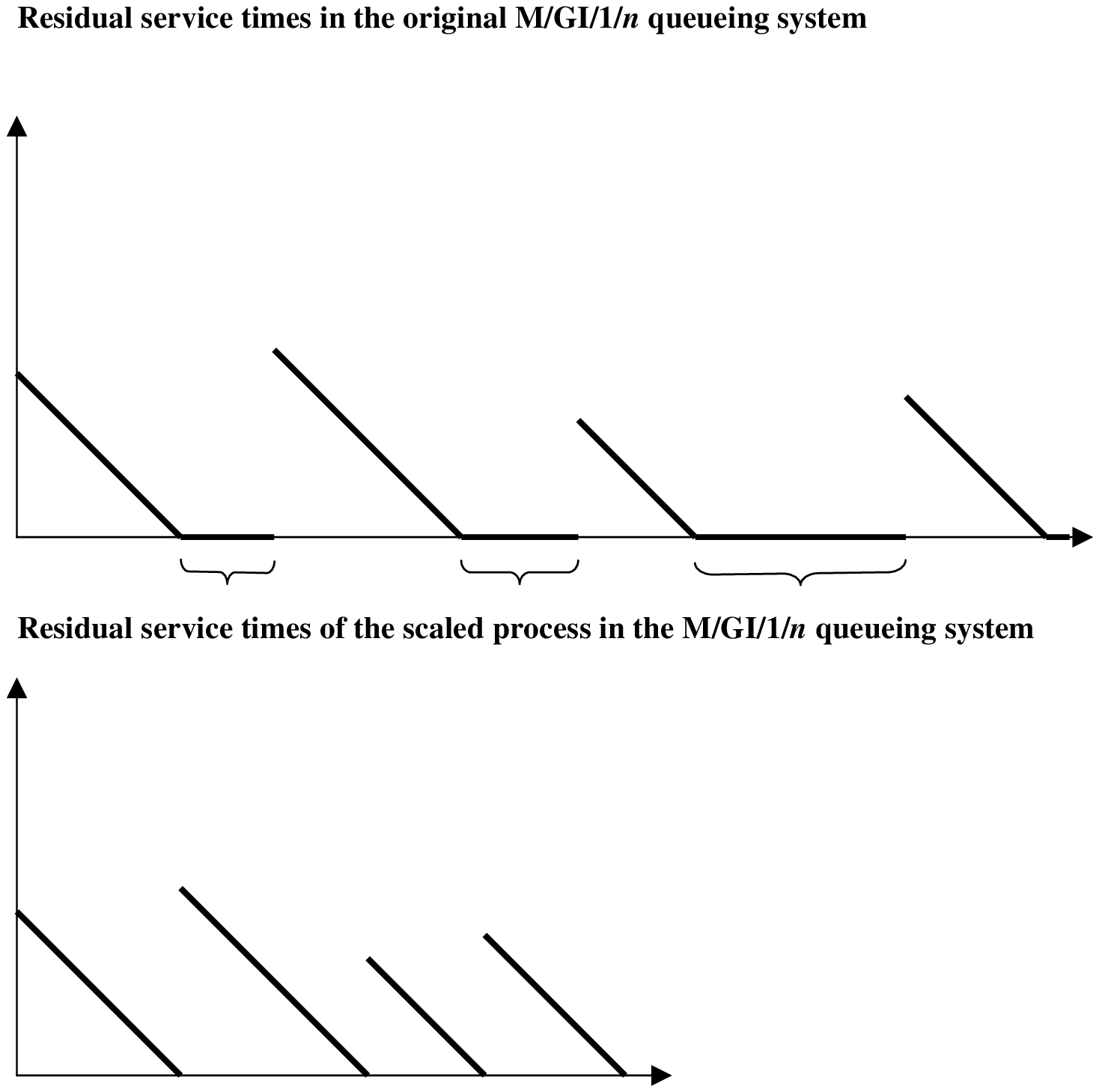}
\caption{Residual service times for the original and scaled
processes of the $M/GI/1/n$ queueing system.}
\end{figure}

In order to establish similar properties in the case $m=2$ let us
first study the properties of 1-stationary processes and explain the construction of \textit{tagged server station} which is substantially used in our construction throughout the paper.

\medskip

\textit{Properties of 1-stationary processes}.
Recall (see Definition \ref{defn1}) that if $\xi(t)$ is a
1-stationary process, then for any $h$ and $t_0$ the probability
distributions of $\xi(t_0)$ and $\xi(t_0+h)$ are the same. The
result remains correct (due to the total probability formula) if
$h$ is replaced by random variable $\vartheta$ with some given
probability distribution, which is assumed to be independent of
the process $\xi$. Namely, we have:
\begin{equation}\label{basic property}
\begin{aligned}
\mathrm{P}\{\xi(t_0+\vartheta)\leq
x\}&=\int_{-\infty}^\infty\mathrm{P}\{\xi(t_0+h)\leq
x\}\mbox{d}\mathrm{P}\{\vartheta\leq h\}\\
&=\mathrm{P}\{\xi(t_0)\leq x\}\int_{-\infty}^\infty
\mbox{d}\mathrm{P}\{\vartheta\leq h\}\\
&=\mathrm{P}\{\xi(t_0)\leq x\}.
\end{aligned}
\end{equation}
That is, $\xi(t_0)$ and $\xi(t_0+\vartheta)$ have the same
distribution.

\smallskip
The above property will be used for the following
\textit{construction of the sequence of 1-stationary processes}
$\xi^{(1)}(t)$, $\xi^{(2)}(t)$, \ldots, having identical
one-dimensional distributions.

Let $\xi^{(0)}(t)=\xi(t)$ be a 1-stationary process, let $t_1$ be
an arbitrary point, and let $\vartheta_1$ be a random variable
with some given probability distribution, which is independent of
the process $\xi^{(0)}(t)$. Let us build a new process
$\xi^{(1)}(t)$ as follows. Put
\begin{equation}\label{trans1}
\xi^{(1)}(t)=\begin{cases} \xi^{(0)}(t), &\mbox{for~all} \
t<t_1,\\
\xi^{(0)}(t+\vartheta_1), &\mbox{for~all} \ t\geq t_1.
\end{cases}
\end{equation}
Since the probability distributions of $\xi(t)$ and
$\xi(t+\vartheta_1)$ are the same for all $t\ge t_1$, then the
processes $\xi(t)$ and $\xi^{(1)}(t)$ have the same
one-dimensional distributions, and $\xi^{(1)}(t)$ is a
1-stationary process as well.

With a new point $t_2$ and a random variable $\vartheta_2$, which
is assumed to be independent of the process $\xi^{(0)}(t)$ and
random variable $\vartheta_1$ (therefore, it is also independent
of the process $\xi^{(1)}(t)$) by the same manner one can build
the new 1-stationary process $\xi^{(2)}(t)$. Specifically,
\begin{equation}\label{trans2}
\xi^{(2)}(t)=\begin{cases} \xi^{(1)}(t), &\mbox{for~all} \
t<t_2,\\
\xi^{(1)}(t+\vartheta_2), &\mbox{for~all} \ t\geq t_2.
\end{cases}
\end{equation}
The new process $\xi^{(2)}(t)$ has the same one-dimensional
distribution as the processes $\xi^{(0)}(t)$ and $\xi^{(1)}(t)$.
The procedure can be infinitely continued, and one can obtain the
infinite family of 1-stationary processes, having the same
one-dimensional distribution.

The points $t_1$, $t_2$,\ldots in the above construction are assumed
to be some fixed (non-random) points. However, the construction also
remains correct in the case of random points $t_0$, $t_1$,\ldots of
Poisson process, since according to the PASTA property \cite{Wolff
1982} the limiting stationary distribution of a 1-stationary process
in a point of a Poisson arrival coincides with the limiting
stationary distribution of the same 1-stationary process in an
arbitrary non-random point. Furthermore, the aforementioned property
of process remains correct when the random points $t_0$,
$t_1$,\ldots are the points of the process which is not necessarily
Poisson but belongs to the special class of processes that contains
Poisson. In this case the property is called ASTA (e.g.
\cite{Melamed}).

\medskip
\textit{1-stationary Poisson process}. Consider an important
particular case when the process $\xi(t)$ is Poisson. Let
$\xi^{(0)}(t)=\xi(t)$. Then the process $\xi^{(1)}(t)$ that obtained
by \eqref{trans1} is no longer Poisson. Its limiting stationary
distribution is the same as that of the original process $\xi(t)$,
but the joint distributions of this process given in different
points $s$ and $t$ distinguish from those of the original process
$\xi(t)$.

The process $\xi^{(1)}(t)$ will be called \textit{1-stationary
Poisson process} or simply \textit{1-Poisson}. Clearly, that the
further processes such as $\xi^{(2)}(t)$, $\xi^{(3)}(t)$, \ldots
that obtained similarly to the procedure in \eqref{trans1},
\eqref{trans2} all are 1-Poisson with the same limiting stationary
distribution. According to the above construction, a 1-Poisson
process is obtained by deleting intervals and merging the ends of an
original Poisson process. Therefore, a sequence of 1-Poisson arrival
time instants is a scaled subsequence of those instants of the
ordinary Poisson arrivals. Hence, for 1-Poisson process the ASTA
property is satisfied, i.e. 1-Poisson arrivals see time averages
exactly as those Poisson arrivals.

\medskip

\textit{Tagged server station}. Consider a stationary queueing
system $M/GI/m/n$, which is referred to as \textit{main server
station}, and in addition to this queueing system introduce another
one containing a server station in order to register specific
arrivals, for example losses or, say, customers waiting their
service in the main system. This server station is called
\textit{tagged server station}. The main idea of introducing tagged
server stations is to decompose the main system as follows. Assume
that along with a Poisson stream of arrivals of customers occupying
servers in the main system, there is another stream of arrivals of
customers in the tagged server system. For instance, the losses in
the main system can be supposed to occupy the tagged server station.
Although the stream of these losses is not Poisson (see e.g.
\cite{Khintchine}, p. 83 or \cite{Kelly - review}, p. 320), it is
shown later that it is 1-Poisson. Therefore, the original system is
decomposed into smaller systems with the same (1-Poisson) type of
input stream. It is worth noting that only one dimensional
distributions of 1-Poisson process are the same for all of them that
generated similarly to the procedure in \eqref{trans1},
\eqref{trans2}. However, the two-dimensional distributions are
distinct in general.

In fact, applications of a tagged server station is wider than that,
and its aim is a proper decomposition of the original system into
the main and tagged systems for further study of the properties of
losses.

Another idea of using tagged server stations is a proper application
of the ASTA property as follows. At the moment of arrival of a
customer in the tagged server station, the stationary
characteristics in the main server station remain the same.
Specifically, the distributions of residual service times in servers
of the main station at the moment of arrival of a customer in the
tagged station coincides with the usual stationary distributions of
these residual service times.

\smallskip
Let us now formulate and prove a property similar to Property
\ref{prop1} for $m=2$. We have the following.

\begin{property}\label{prop2}
For the $M/GI/2/0$ queueing system we have:
\begin{equation}\label{fact1}
\mathrm{P}\{\widehat{{\bf y}}_{1,0}(t)\in B_1\}=\mathrm{P}\{{\bf
x}_1(t)\in B_1\},
\end{equation}
\end{property}

\begin{proof}
In order to simplify the explanation in this case, let us consider
two auxiliary stationary one-dimensional processes $\zeta_{1,0}(t)$
and $\zeta_{2,0}(t)$. The first process describes a residual service
time in the first server, and the second one describes a residual
service time in the second server. If the $i$th server ($i=1,2$) is
free in time $t$, then we set $\zeta_{i,0}(t):=0$.

Our further convention is that the first server is a \textit{tagged}
server. We assume that if at the moment of arrival of a customer
both of the servers are free, he/she occupies the first server.
Clearly that this assumption is not a loss of generality. For
instance, if we assume that both of the servers are equivalent and
can be occupied with the equal probability $\frac{1}{2}$, then an
occupied server (let it be the first) can be called tagged. In
another busy period start an arriving customer can occupy the second
server. It this case, nothing is changed if the servers will be
renumbered, and the occupied server will be numbered as first and
called tagged.

%
Our main idea is a decomposition of the stationary $M/GI/2/0$
queueing system into two systems and study the properties of
stationary (1-stationary) processes $\zeta_{1,0}(t)$ and
$\zeta_{2,0}(t)$. The arrival stream to the tagged system is
Poisson, so the first system is $M/GI/1/0$, while the second one is
denoted $\bullet/GI/1/0$, where $\bullet$ in the first place of the
notation stands for the input process in the second system, which is
the output (loss) stream in the first one. Clearly, that an arriving
customer is arranged to the second queueing system if and only if at
the moment of his/her arrival the tagged system is occupied.
Therefore, let us delete all the intervals when the tagged system is
empty and merge the ends. In this case, the tagged system becomes an
ordinary renewal process, and the stream of arrivals to the second
queueing system becomes 1-Poisson rather then Poisson (because after
deleting intervals and merging the ends in the new time scale the
original Poisson process is transformed into 1-Poisson). Therefore
the second system now can be re-denoted by $\widetilde{M}/GI/1/0$,
where $\widetilde{M}$ in the first place of the notation stands for
1-Poisson input and replaces the initially written symbol $\bullet$.

Thus, the $M/GI/2/0$ queueing system is decomposed into the
$M/GI/1/0$ and $\widetilde{M}/GI/1/0$ queueing systems. Clearly,
that without loss of generality one can assume that the original
arrival stream is 1-Poisson rather than Poisson, i.e. the original
queueing system is $\widetilde{M}/GI/2/0$, and it is decomposed into
two $\widetilde{M}/GI/1/0$ queueing systems. The last note is
important for the further extension of the result for the systems
$M/GI/m/0$ (or generally $\widetilde{M}/GI/m/0$) having $m>2$
servers.




Let $\tau$ be the time moment when an arriving customer occupies the
tagged server station.
According
to the ASTA property,
\begin{equation}\label{add-prop1}
\mathrm{P}\{\zeta_{2,0}(\tau)\leq x\}=\mathrm{P}\{\zeta_{2,0}(t)\leq
x\},
\end{equation}
where $t$ is an arbitrary fixed point, and the probability
distribution function of $\zeta_{2,0}(t)$ in this point coincides
with the distribution of residual service time in specific
$\widetilde{M}/GI/1/0$ system with some specific value of parameter
of 1-Poisson process, which is not important here. On the other
hand, the process $\zeta_{2,0}(t)$ is stationary and Markov.
Therefore from \eqref{add-prop1} for any $h>0$ we have
\begin{equation}\label{add-prop3}
\begin{aligned}
\mathrm{P}\{\zeta_{2,0}(\tau+h)\leq
x\}&=\mathrm{P}\{\zeta_{2,0}(t+h)\leq
x\}=\mathrm{P}\{\zeta_{2,0}(t)\leq x\}.
\end{aligned}
\end{equation}

Let $\chi$ denotes the service time of the customer, who arrives at
the time moment $\tau$ occupying the tagged server station. Our
challenge is to prove that
\begin{equation}\label{add-prop7}
\mathrm{P}\{\zeta_{2,0}(\tau+\chi)\leq
x\}=\mathrm{P}\{\zeta_{2,0}(t)\leq x|\zeta_{1,0}(t)>0\}.
\end{equation}
Instead of the original processes $\zeta_{i,0}(t)$, $i=1,2$,
consider another processes $\widetilde{\zeta}_{i,0}(t)$, which are
obtained by deleting the intervals where the tagged server is free,
and merging the ends. Then, $\widetilde{\zeta}_{1,0}(t)$ is a
renewal process, and the 1-stationary process
$\widetilde{\zeta}_{2,0}(t)$ and the random variable $\chi$ (the
length of a service time in the tagged server that starts at moment
$\tau$) are independent.
%
%
%
%
%
Hence, for
any event $\{\chi=h\}$ according to the properties of 1-stationary
processes we have
\begin{equation}\label{add-prop4}
\begin{aligned}
\mathrm{P}\{\widetilde{\zeta}_{2,0}(\tau+\chi)\leq
x|\chi=h\}&=\mathrm{P}\{\widetilde{\zeta}_{2,0}(\tau)\leq x\},
\end{aligned}
\end{equation}
and, due to the total probability formula from \eqref{add-prop4} we
have
\begin{equation}\label{add-prop2}
\mathrm{P}\{\widetilde{\zeta}_{2,0}(\tau+\chi)\leq
x\}=\mathrm{P}\{\widetilde{\zeta}_{2,0}(\tau)\leq x\}.
\end{equation}
The only difference between \eqref{add-prop2} and the basic property
\eqref{basic property} is that the time moment $\tau$ is random,
while $t_0$ is not. However keeping in mind \eqref{add-prop3}, this
modified equation \eqref{add-prop2} follows by the same derivation
as in \eqref{basic property}.

Hence, from \eqref{add-prop2},
\begin{equation}\label{add-prop5}
\mathrm{P}\{\widetilde{\zeta}_{2,0}(\tau+\chi)\leq
x\}=\mathrm{P}\{\widetilde{\zeta}_{2,0}(\tau)\leq
x\}=\mathrm{P}\{\zeta_{2,0}(t)\leq x|\zeta_{1,0}(t)>0\},
\end{equation}
and since $\mathrm{P}\{\widetilde{\zeta}_{2,0}(\tau+\chi)\leq
x\}=\mathrm{P}\{{\zeta}_{2,0}(\tau+\chi)\leq x\}$ we finally arrive
at \eqref{add-prop7}.

As well, noticing that
$$
\mathrm{P}\{\widetilde{\zeta}_{2,0}(\tau)\leq
x\}=\mathrm{P}\{{\zeta}_{2,0}(\tau)\leq x\},
$$
from \eqref{add-prop2} and \eqref{add-prop1} we also have
\begin{equation}\label{ASTA-basic}
\mathrm{P}\{{\zeta}_{2,0}(\tau+\chi)\leq
x\}=\mathrm{P}\{{\zeta}_{2,0}(\tau)\leq
x\}=\mathrm{P}\{\zeta_{2,0}(t)\leq x\}.
\end{equation}


Similarly to \eqref{ASTA-basic} one can prove
\begin{equation}\label{ASTA-basic2}
\mathrm{P}\{{\zeta}_{1,0}(\tau+\chi)\leq
x\}=\mathrm{P}\{{\zeta}_{1,0}(\tau)\leq
x\}=\mathrm{P}\{\zeta_{1,0}(t)\leq x\},
\end{equation}
where $\tau$ is the moment of arrival of a customer, who at this
moment $\tau$ occupies the second server, and $\chi$ denotes his/her
service time. Relations \eqref{ASTA-basic2} can be proved with the
aid of the same construction of deleting intervals and merging the
ends, but now in the second server. So, combining \eqref{ASTA-basic}
and \eqref{ASTA-basic2} we arrive at the following fact. \textit{In
any arrival or service completion time instant in one server},
\textit{the residual service time in another server has the same
stationary distribution}.

This fact is used in the constructions below.

Now consider the stationary $M/GI/2/0$ queueing system, in which
both servers are equivalent in the sense that if at the moment of
arrival of a customer both servers are free, then a customer can
occupy each of servers with the equal probability $\frac{1}{2}$. In
this case the both of the processes $\zeta_{1,0}(t)$ and
$\zeta_{2,0}(t)$ have the same distribution.

Let us delete the time intervals where the both servers are
simultaneously free, and merge the corresponding ends (see Figure
3). The new processes are denoted by $\widetilde\zeta_{1,0}(t)$ and
$\widetilde\zeta_{2,0}(t)$, and both of them have the same
equivalent distribution. (We use the same notation as in the
construction above believing that it is not confusing for readers.)
This two-dimensional 1-stationary process characterizes the system
where in any time $t$ at least one of two servers is busy. Consider
the event of arrival of a customer in a stationary system at the
moment when only one of two servers is busy. Let $\tau^*$ be the
moment of this arrival, and let $\tau^{**}$ denote the moment of the
first service completion in one of two servers following after the
moment $\tau^*$. Then, at the endpoint $\tau^{**}$ of the interval
$[\tau^*, \tau^{**})$ the distribution of the residual service time
will be the same as that at the moment $\tau^*$ (due to the
established fact that at the end of a service completion in one
server, the distribution of a residual service time in another
server must coincide with the stationary distribution of a residual
service time and due to the fact that both servers are equivalent.)

The additional details here are as follows. There can be different
events associated with the points $\tau^*$ and $\tau^{**}$. For
example, at the moment $\tau^*$ an arriving customer can be accepted
by one of the servers, while the service completion at the moment
$\tau^{**}$ can be either in the same server of in another server.
If time moments $\tau^*$ and $\tau^{**}$ are associated with the
same server (for example, the moment of service start and service
completion in the first server) then we speak about residual service
times in another server (in this example - the second server). If
time moments $\tau^*$ and $\tau^{**}$ are associated with different
servers (say, $\tau^*$ is the service start in the first server, but
$\tau^{**}$ is the service completion in the second one), then we
speak about residual service times in different servers (in this
specific case we speak about residual service time in the second
server at the time moment $\tau^*$ and residual service time in the
first server at the time moment $\tau^{**}$). However, according to
the earlier result, it does not matter which specific event of these
mentioned occurs. The only fact, that the stationary distribution of
a residual service time in a given server must be the same for all
time moments of arrival and service completion occurring in another
server and vice versa, is used.

Deleting the interval [$\tau^*$, $\tau^{**}$) and merging the ends
$\tau^*$ and $\tau^{**}$ (see Figure 4) we obtain the following
structure of the 1-stationary process $\widehat{{\bf y}}_{1,0}(t)$.

In the points where idle intervals are deleted and the ends are
merged we have renewal points: one of periods is finished and
another is started. In the other points where the intervals of type
[$\tau^*$, $\tau^{**}$) are deleted and their ends are merged we
have the points of `interrupted' renewal processes. In this
`interrupted' renewal process the point $\tau^*$ is a point of
1-Poisson arrival, and, according to ASTA, the distribution in this
point in the server that continue to serve a customer coincides with
the stationary distribution of the residual service time. In the
other point $\tau^{**}$, which is the point of a service completion,
the distribution in this point in the server that continue to serve
a customer coincides with the stationary distribution of a residual
service time as well. Therefore, in the point of the interruption
(which is a point of discontinuity) the residual service time
distribution coincides with the stationary distribution of a
residual service time, i.e. with the distribution of ${\bf x}_1(t)$.
(Notice, that the intervals of type [$\tau^*, \tau^{**}$) are an
analogue of the intervals [$s_{1,k}, t_{1,k}$) considered in the
Markovian case in Section \ref{Markovian}.)

By amalgamating the residual service times of the first and second
servers given in the lower graph in Figure 4, one can built a
typical one-dimensional 1-stationary process $\widehat{{\bf
y}}_{1,0}(t)$, the limiting stationary distribution of which
coincides with that of ${\bf x}_1(t)$. (see Figure 5).

Therefore the processes $\widehat{{\bf y}}_{1,0}(t)$ and ${\bf
x}_1(t)$ have the identical one-dimensional distribution, and
relation \eqref{fact1} follows.
\end{proof}

\begin{figure}\label{fig2}
\includegraphics[width=15cm,height=20cm]{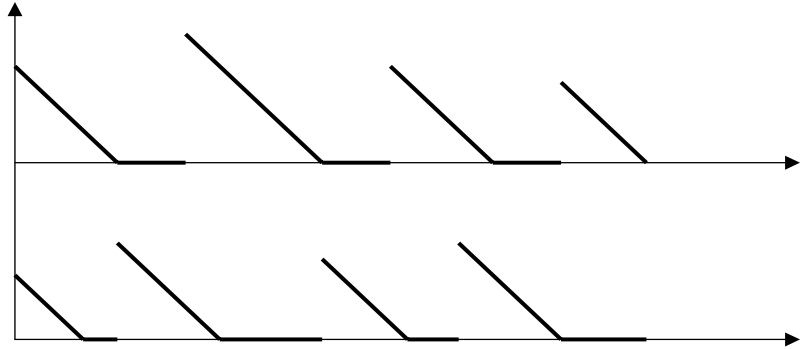}
\caption{Residual service times for the original and scaled
processes of the $M/GI/2/0$ queueing system after deleting
intervals where both of the servers are free, and merging the
ends.}
\end{figure}
\begin{figure}\label{fig3}
\includegraphics[width=15cm,height=20cm]{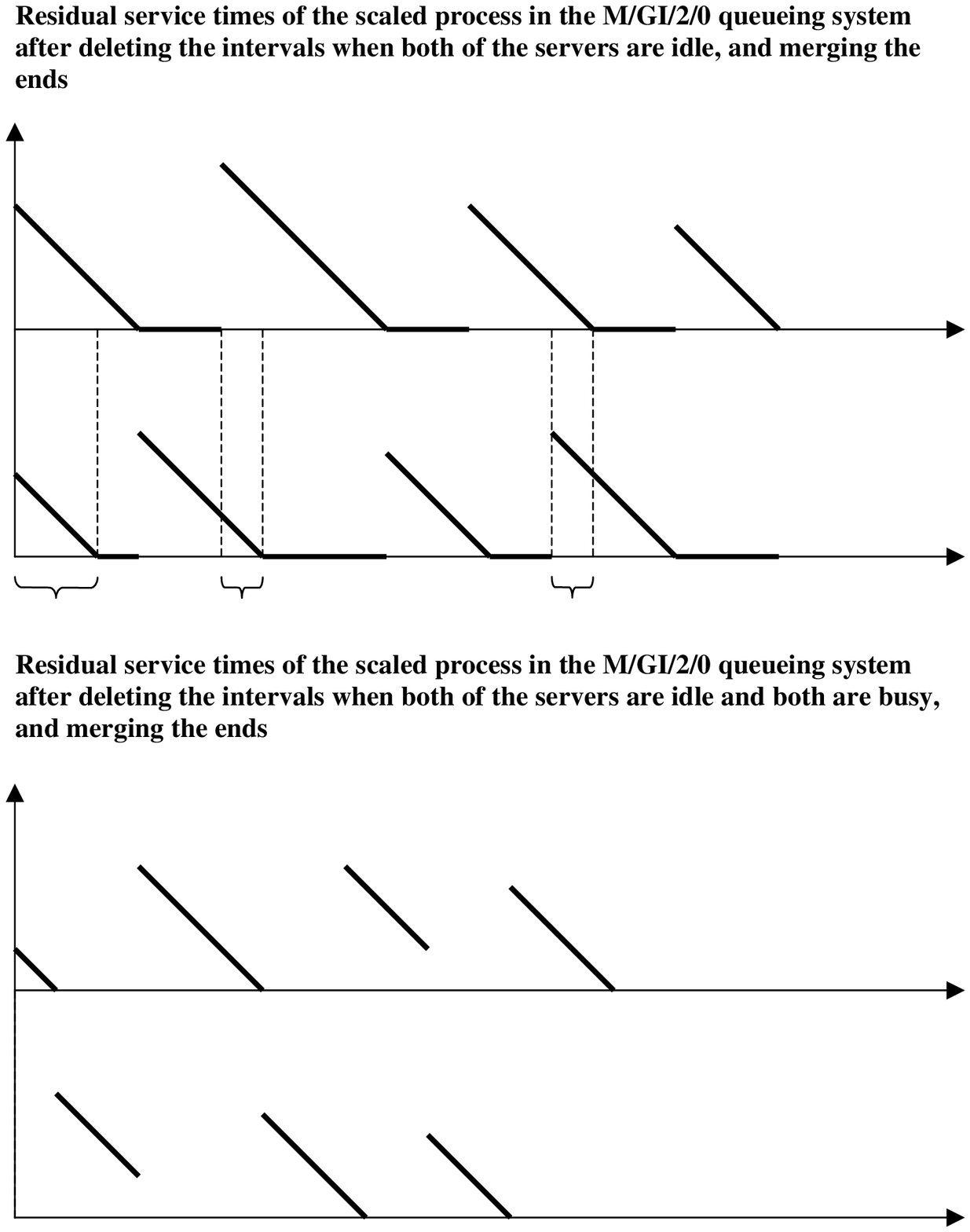}
\caption{Residual service times for the original and scaled
processes of the $M/GI/2/0$ queueing system after deleting
intervals where both of the servers are free and both are busy,
and merging the ends.}
\end{figure}
\begin{figure}\label{figAmalg}
\includegraphics[width=15cm,height=20cm]{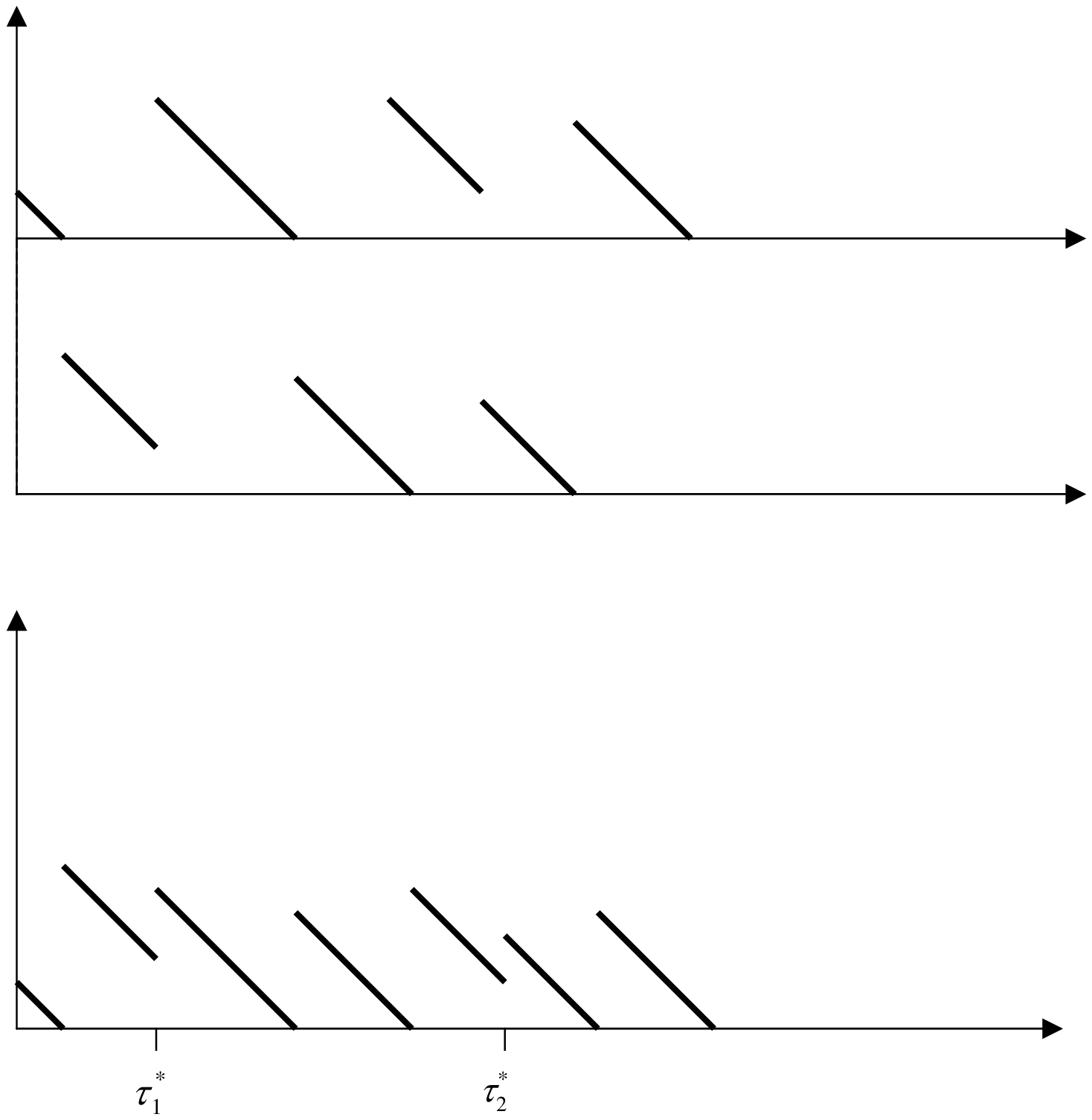}
\caption{A typical 1-stationary process of residual service times
is obtained by amalgamating residual service times of the first
and second servers. $\tau_1^*$ and $\tau_2^*$ are the points where
the intervals of type [$\tau^*, \tau^{**}$) are deleted, and the
ends are merged}
\end{figure}

Let us develop Property \ref{prop2} to the case $m=3$ and then to
the case of an arbitrary $m>1$ for the $M/GI/m/0$ queueing
systems. Namely, we have the following.

\begin{property}\label{prop3} For the $M/GI/m/0$ queueing system we have:
\begin{equation}\label{fact1-ext}
\mathrm{P}\{\widehat{{\bf y}}_{m-1,0}(t)\in
B_{m-1}\}=\mathrm{P}\{{\bf x}_{m-1}(t)\in B_{m-1}\},
\end{equation}
where $B_{m-1}$ is an arbitrary Borel set of $\mathbb{R}^{m-1}$.
\end{property}

\begin{proof} The proof will be concentrated in the case $m=3$ for
the 1-stationary process $\widehat{{\bf y}}_{2,0}(t)$, which is
associated with the paths of the $M/GI/3/0$ queueing system where
only two servers are busy. Then the result will be concluded for an arbitrary $m\geq2$ by induction.

Prior studying this case, we first study the specific case of the
$M/GI/2/0$ queueing system by considering the paths when the both
servers are busy. Then using the arguments of the proof of Property
\ref{prop2} enables us to extend that specific result related to the
$M/GI/2/0$ queueing systems to the 1-stationary process
$\widehat{{\bf y}}_{2,0}(t)$ of the $M/GI/3/0$ queueing system.

As in the proof of Property \ref{prop2} in the specific case of the
$M/GI/2/0$ queueing system considered here, we will study the
stationary one-dimensional processes $\zeta_{1,0}(t)$ and
$\zeta_{2,0}(t)$. However the idea of the present proof generally
differs from that of the proof of Property \ref{prop2}. Here we do
not call the first (or second) server a tagged server station to use
decomposition. We simply use the fact established in the proof of
Property \ref{prop2} that at the moment of arrival or service
completion of a customer in one server, the distribution of a
residual service time in another server will coincide with the
stationary distribution of a residual service time in this server.
(The same idea has been used in the proof of Property \ref{prop2}.)

The present proof
explicitly uses the fact that the class of 1-stationary processes
is algebraically closed with respect to the operations of deleting
intervals and merging the ends, which was mentioned before.

Let us delete the idle intervals of the process $\zeta_{1,0}(t)$
and merge the ends. Then we get a stationary renewal process as in
the above case $m=1$ (Property \ref{prop1}).

After deleting the same time intervals in the second stationary
process $\zeta_{2,0}(t)$ and merging the ends, the process will be
transformed as follows. Let $t^*$ be a moment of 1-Poisson arrival
when a customer occupies the first server. (Recall that owing to the
known properties of 1-Poisson process, the stream of arrival to each
of $i$ servers ($i=1,2$) is 1-Poisson.) Then, according to the ASTA
property, $\zeta_{2,0}(t^*)=\zeta_{2,0}(t)$ in distribution.
Therefore after deleting all of the idle intervals of the second
server and merging the ends, after the first time scaling (i.e.
removing corresponding time intervals, see Figure 6) instead of the
initial 1-stationary process $\zeta_{2,0}(t)$ we obtain the new
1-stationary process with the equivalent one-dimensional
distribution. This process is denoted by $\widetilde\zeta_{2,0}(t)$.

\begin{figure}\label{fig4}
\includegraphics[width=15cm,height=20cm]{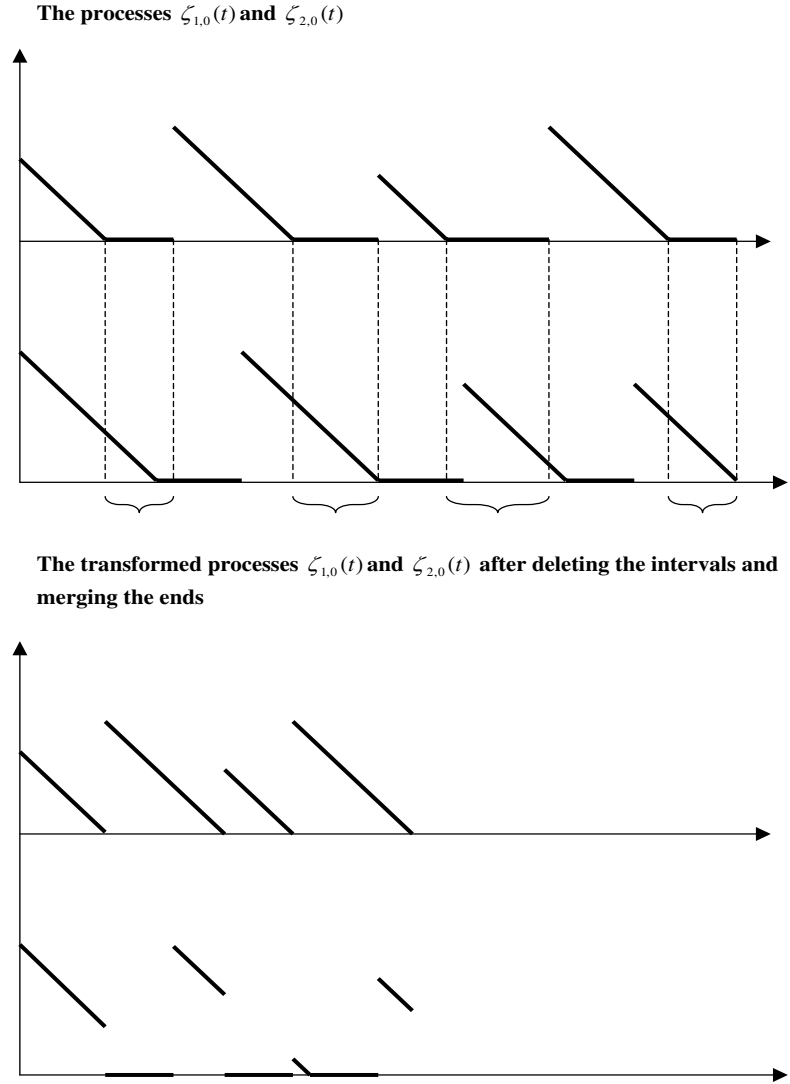}
\caption{The dynamic of time scaling for a queueing system with
two servers after deleting the idle intervals in the first server
and merging the ends}
\end{figure}

Notice, that the process $\widetilde\zeta_{2,0}(t)$ is obtained
from the process $\zeta_{2,0}(t)$ by constructing a sequence of
1-stationary processes described above.

Then we have the two-dimensional process the first component of
which is ${\bf x}_1(t)$ and the second one is
$\widetilde\zeta_{2,0}(t)$. For our convenience this first
component is provided with upper index, and the two-dimensional
vector looks now as $\left\{{\bf x}_1^{(1)}(t),
\widetilde\zeta_{2,0}(t)\right\}$.

Let us repeat the above procedure, deleting the remaining idle
intervals of the second server and merging the ends. We get the
1-stationary process being equivalent in the distribution to the
stationary renewal process ${\bf x}_1(t)$, which is denoted now
${\bf x}_1^{(2)}(t)$.

Upon this (final) time scaling the first process ${\bf
x}_1^{(1)}(t)$ is transformed as follows. Let $t^{**}$ be a random
point of 1-Poisson arrival when the second server is occupied.
Applying the ASTA property once again, for the first component of
the process we obtain that ${\bf x}_1^{(1)}(t^{**})$ coincides in
one-dimensional distribution with ${\bf x}_1^{(1)}(t)$. Thus, after
deleting the entire idle intervals and merging the ends, we finally
obtain the two-dimensional process $\left\{{\bf x}_1^{(1)}(t), {\bf
x}_1^{(2)}(t)\right\}$ each component of which has the same
one-dimensional distribution as this of the process ${\bf x}_1(t)$.
The dynamic of this time scaling is shown in Figure 7.

\begin{figure}\label{fig5}
\includegraphics[width=15cm,height=20cm]{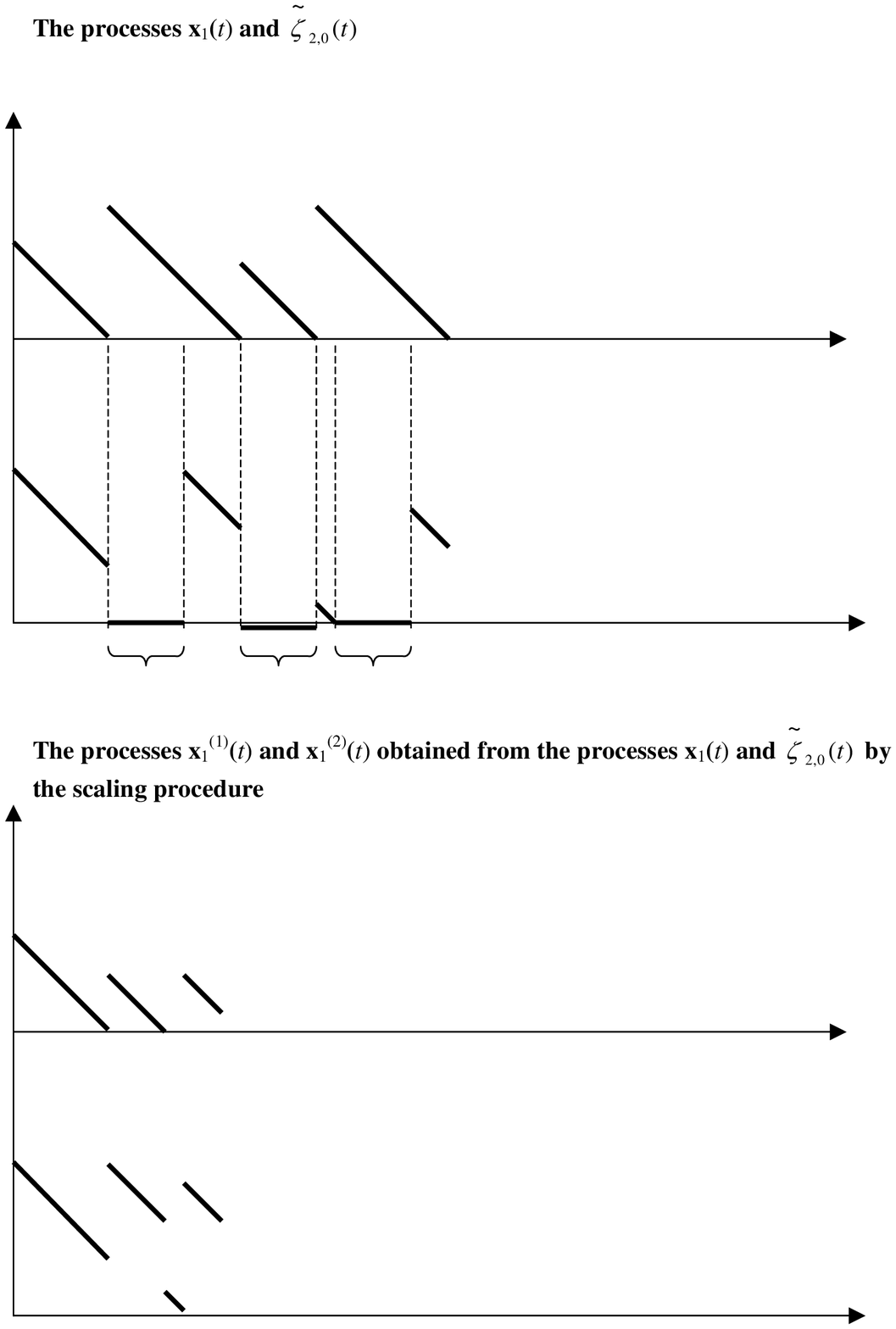}
\caption{The dynamic of time scaling for a queueing system with
two servers after deleting the idle intervals in the second server
and merging the ends}
\end{figure}

For our further purpose, the independence of the processes ${\bf
x}_1^{(1)}(t)$ and ${\bf x}_1^{(2)}(t)$ is needed. The constructions
in this paper enables us to prove this independence. However, the
independence of ${\bf x}_1^{(1)}(t)$ and ${\bf x}_1^{(2)}(t)$
follows automatically from the known results by Tak\'acs
\cite{Takacs 1969} and the easiest way is to follow a result of that
paper. Namely, it follows from formulae (6) and (7) on page 72, that
the joint conditional stationary distribution of residual service
times given that $k$ servers are busy coincides with the stationary
distribution of ${\bf x}_k(t)$, which in turn is the product of the
stationary distributions of ${\bf x}_1(t)$. In particular,
\begin{equation}\label{3.9}
{\rm P}\left\{{\bf x}_1^{(1)}(t)\leq x_1, \ {\bf x}_1^{(2)}(t)\leq
x_2\right\}={\rm P}\left\{{\bf x}_1(t)\leq x_1\right\}{\rm
P}\left\{{\bf x}_1(t)\leq x_2\right\}.
\end{equation}

Now, using the arguments of the proof of Property \ref{prop2} one
can easily extend the result obtained now for $M/GI/2/0$ queueing
system to the $M/GI/3/0$ queueing system, and thus prove
\eqref{fact1-ext} for the $M/GI/3/0$ queueing system.

Similarly to the proof of Property \ref{prop2}, let us introduce the
processes $\zeta_{1,0}(t)$, $\zeta_{2,0}(t)$ and $\zeta_{3,0}(t)$ of
the residual service times in the first, second and third servers
correspondingly. These processes all are assumed to have the same
stationary distribution of residual times, which respects to the
scheme where an arriving customer occupies one of available free
servers with equal probability. Let us call a server of the
${M}/GI/3/0$ queueing system that occupied at the moment of busy
period start a tagged server station. So, we decompose the original
system into the $\widetilde{M}/GI/2/0$ and tagged queueing system
$\widetilde{M}/GI/1/0$. However, it is shown above that
$\widetilde{M}/GI/2/0$ can be decomposed into two
$\widetilde{M}/GI/1/0$ queueing systems, where after the procedure
of deleting idle intervals and merging the ends we obtain the
process having the same stationary distribution as that of the
process $\mathbf{x}_2(t)$. This stationary distribution remains the
same in all random points of arrivals and service completions in the
tagged service station. So, after deleting intervals and merging the
ends in the tagged service station, in a new time scaling we arrive
at the same stationary distribution as that of the process
$\mathbf{x}_3(t)$. So, the result for $m=3$ follows.

This induction becomes clear for an arbitrary $m\geq2$ as well,
where the original $\widetilde{M}/GI/m/0$ system can be decomposed
into the $\widetilde{M}/GI/m-1/0$ queueing system and a tagged
server station $\widetilde{M}/GI/1/0$.
\end{proof}

\smallskip
Now we will establish a connection between the processes
$\widehat{{\bf y}}_{m-1,n}(t)$, $\widehat{{\bf y}}_{m-1,0}(t)$ and
${\bf x}_{m-1}(t)$. We start from the case $m=2$.

\begin{property}\label{prop4}
Under the assumption that the probability distribution function
$G(x)$ belongs to the class NBU (NWU) we have
\begin{equation}\label{fact2}
\mathrm{P}\{\widehat{{\bf y}}_{1,n}(t)\leq x\}\leq \ (\mbox{resp.}
\geq)\ \mathrm{P}\{{\bf x}_1(t)\leq x\}.
\end{equation}
\end{property}
\begin{proof}
Along with the 1-stationary processes $\widehat{\bf y}_{1,n}(t)$,
let us introduce another 1-stationary processes $\widehat{\bf
Y}_{2,n}(t)$. This last process is related to the same
$\widetilde{M}/GI/2/n$ queueing system as the process $\widehat{\bf
y}_{1,n}(t)$, and is obtained by deleting time intervals when there
are more than two or less than two customers in the system, and
merging the ends.

Using the same arguments as in the proof of Property \ref{prop3} one
can prove that the components of this process are generated by
independent 1-stationary processes and having the same distribution.
Indeed, involving as earlier in the proof of Property \ref{prop2}
the processes $\zeta_1(t)$ and $\zeta_2(t)$ having the same
distribution, one can delete intervals where the system is empty and
merge the ends. Apparently, the new processes $\widetilde\zeta_1(t)$
and $\widetilde\zeta_2(t)$ obtained after this procedure have the
same stationary distribution. (However, it is shown later that the
limiting stationary distribution of these one-dimensional processes
differs from such the distribution obtained after the similar
procedure for the $\widetilde{M}/GI/2/0$ system, and, therefore, its
one-dimensional distribution distinguishes from that of the process
${\bf x}_1(t)$.)

Let us go back to the initial process ${\bf y}_{2,n}(t)$,
$n\geq1$, to delete the time intervals where the both servers are
free and merge the corresponding ends. We also remove the last
component corresponding to the queue-length $Q_{2,n}(t)$. (The
exact value of the queue-length $Q_{2,n}(t)$ is irrelevant here
and is not used in our analysis.) In the new time scale we obtain
the two-dimensional process $\widetilde{{\bf y}}_{2,n}(t)$.

Similarly to the proof of relation \eqref{fact1} we have time
moments $\tau^{*}$ and $\tau^{**}$. The first of them is a moment
of arrival of a customer at the system with one busy and one free
server, and the second one is the following after $\tau^*$ moment
of service completion of a customer when there remain one busy
server only. The time interval [$\tau^*$, $\tau^{**}$) is an
\textit{orbital busy period}. (The concept of orbital busy period
is defined in Section \ref{Markovian} for Markovian systems. For
$M/GI/m/n$ queueing systems this concept is the same.) It can
contain queueing customers waiting for their service. Let
$t^{\mathrm{begin}}$ be a moment of arrival of a customer during
the orbital busy period [$\tau^*$, $\tau^{**}$) who occupies a
waiting place, and let $t^{\mathrm{end}}$ be the following after
$t^{\mathrm{begin}}$ moment of time when after the service
completion the queue space becomes empty again. A period of time
[$t^{\mathrm{begin}}$, $t^{\mathrm{end}}$) is called
\textit{queueing period}. (Note that for the M/GI/$m$/$n$ queueing
system, the intervals of type [$\tau^*$, $\tau^{**}$) are an
analogue of the intervals of type [$s_{m-1,k}, t_{m-1,k}$) in the
Markovian queueing system $M/M/m/n$, and the intervals of type
[$t^{\mathrm{begin}}, t^{\mathrm{end}}$) are an analogue of the
intervals of type [$s_{m,k}, t_{m,k}$) in that Markovian queueing
system $M/M/m/n$.)

All customers of \textit{queueing periods}, i.e. those arrived during orbital busy period can be considered as customers arriving in a tagged server station. At the moment of $t^{\mathrm{begin}}$, which is an instant of a Poisson arrival, the two-dimensional distribution of the random vector $\widetilde{{\bf
y}}_{2,n}(t^{\mathrm{begin}})$ coincides with the stationary distribution of the random vector $\widehat{{\bf Y}}_{2,n}(t)$. However, in the point $t^{\mathrm{end}}$, the probability
distribution of $\widetilde{{\bf y}}_{2,n}(t^{\mathrm{end}})$ is
different from the stationary distribution of $\widehat{{\bf Y}}_{2,n}(t)$, because this specific time instant $t^{\mathrm{end}}$ coincides with a service beginning in one of servers of the main system. Therefore, deleting the interval [$t^{\mathrm{begin}}, t^{\mathrm{end}}$) and merging the end leads to the change of the distribution.

More specifically, at the time instant $t^{\mathrm{end}}$ one of the
components of the vector $\widehat{{\bf Y}}_{2,n}(t)$, say the first
one, is a random variable having the probability distribution
$G(x)$. Then, another component, i.e. the second one, because of the
aforementioned properties of 1-stationary processes, coincides in
distribution with $\widetilde\zeta_{1,n}$ (or
$\widetilde\zeta_{2,n}$), which is a component of the stationary
process $\widehat{{\bf Y}}_{2,n}(t)$. Indeed, let customers arriving
in a busy system and waiting in the queue be assigned to the tagged
server station. At the moment of 1-Poisson arrival
$t^{\mathrm{begin}}$ of a customer in the tagged server station, the
two-dimensional Markov process associated with the main queueing
system has the same distribution as the vector $\widehat{{\bf
Y}}_{2,n}(t)$, i.e. two-dimensional distribution coinciding with the
joint distribution of $\widetilde\zeta_{1,n}$ and
$\widetilde\zeta_{2,n}$. Then, at the moment of the service
completion $t^{\mathrm{end}}$, which coincides with the moment of
service completion in one of two servers, the probability
distribution function of the residual service time in another
server, where the service is being continued, coincides with the
distribution of a component of the vector $\widehat{{\bf
Y}}_{2,n}(t)$, i.e. with the distribution of
$\widetilde\zeta_{1,n}$.

If the probability distribution function $G(x)$ belongs to the class NBU,
then the 1-stationary process $\widetilde{{\bf y}}_{2,n}(t)$ satisfies
the property $\widetilde{{\bf
y}}_{2,n}(t^{\mathrm{begin}})\leq_{st}\widetilde{{\bf
y}}_{2,n}(t^{\mathrm{end}})$. If $G(x)$ belongs to the class NWU, then the opposite inequality holds:
$\widetilde{{\bf
y}}_{2,n}(t^{\mathrm{end}})\leq_{st}\widetilde{{\bf
y}}_{2,n}(t^{\mathrm{begin}})$.
(The stochastic inequality between
vectors means the stochastic inequality between their
corresponding components.)

The above stochastic inequalities are between random values of the
process $\widetilde{{\bf y}}_{2,n}(t)$ in the different time
instants $t^{\mathrm{begin}}$ and $t^{\mathrm{end}}$. Our further
task is to compare two different processes $\widetilde{{\bf
y}}_{2,n}(t)$ and $\widetilde{{\bf y}}_{2,0}(t)$. The first of these
processes is associated with the $M/GI/m/n$ queueing system, while
the second one is associated with the $M/GI/m/0$ queueing system.
The idea of comparison is very simple. Suppose that both queueing
system are started at zero, i.e. consider the paths of these system
when the both of them are not in steady state, and compare the
Markov processes associated with these system. For the
non-stationary processes we will use the same notation
$\widetilde{{\bf y}}_{2,n}(t)$ and $\widetilde{{\bf y}}_{2,0}(t)$
understanding that it is spoken about usual (not stationary) Markov
processes. The notation for time moments such as
$t^{\mathrm{begin}}$ and $t^{\mathrm{end}}$ is now associated with
these usual (i.e. non-stationary) processes as well. We will
consider the Markov processes associated with $M/GI/2/n$ and
$M/GI/2/0$ queueing systems on the same probability space. In the
time interval $[0, t^{\mathrm{begin}})$ the paths of the Markov
processes $\widetilde{{\bf y}}_{2,n}(t)$ and $\widetilde{{\bf
y}}_{2,0}(t)$ coincide ($n\neq 0$). However, after deleting the
interval $[t^{\mathrm{begin}}, t^{\mathrm{end}})$ and merging the
ends, then in the end point $t^{\mathrm{begin}}$ the values of the
processes $\widetilde{{\bf y}}_{2,n}(t)$ and $\widetilde{{\bf
y}}_{2,0}(t)$ will be different. Indeed, in the case of the process
$\widetilde{{\bf y}}_{2,0}(t)$, which is associated with the
$M/GI/m/0$ queueing system, $t^{\mathrm{begin}}$ and
$t^{\mathrm{end}}$ is the same point, and the value of Markov
processes will be the same after replacing the points
$t^{\mathrm{begin}}$ with $t^{\mathrm{end}}$. However, in the case
of the process $\widetilde{{\bf y}}_{2,n}(t)$ associated with
$M/GI/m/n$ queueing system, the values in these points will be
different with probability 1, and consequently, because of the
inequality $\widetilde{{\bf
y}}_{2,n}(t^{\mathrm{begin}})\leq_{st}\widetilde{{\bf
y}}_{2,n}(t^{\mathrm{end}})$ we have $\widetilde{{\bf
y}}_{2,0}(t^{\mathrm{begin}})\leq_{st}\widetilde{{\bf
y}}_{2,n}(t^{\mathrm{begin}})$ in the case when $G(x)$ belongs to
the class NBU. If $G(x)$ belongs to the class NWU, we have the
opposite inequality: $\widetilde{{\bf
y}}_{2,n}(t^{\mathrm{begin}})\leq_{st}\widetilde{{\bf
y}}_{2,0}(t^{\mathrm{begin}})$.

Therefore, after deleting all the intervals of the type
$[t^{\mathrm{begin}}, t^{\mathrm{end}})$ from the original Markov
process we obtain new Markov process, and in the case when $G(x)$
belongs either to the class NBU or to the class NWU one can apply
the theorem of Kalmykov \cite{Kalmykov 1962} (see also \cite{Keilson
1979}) to compare these two Markov processes. In the case where
$G(x)$ belongs to the class NBU, all the path of the Markov process,
associated with $M/GI/2/n$ is not smaller (in stochastic sense) than
that path of the Markov process, associated with $M/GI/2/0$. If
$G(x)$ belongs to the class NWU, then the opposite stochastic
inequality holds between two different Markov processes. Apparently,
the same stochastic inequalities remain correct if we speak about
stationary Markov processes. Nothing is changed if we let $t$ to
increase to infinity and arrive at stationary distributions. So,
under the assumption that $G(x)$ belongs to the class NBU, for the
stationary processes we obtain $\widehat{{\bf
y}}_{1,0}(t)\leq_{st}\widehat{{\bf y}}_{1,n}(t)$. In other words,
due to the fact that $\widehat{{\bf y}}_{1,0}(t)=_{st}{\bf x}_1(t)$,
we obtain that ${\bf x}_1(t)\leq_{st}\widehat{{\bf y}}_{1,n}(t)$. In
the case where $G(x)$ belongs to the class NWU, the opposite
inequality holds.

The arguments of the proof given for $m=2$ remain correct for an arbitrary $m\geq2$. The proof given by induction uses decomposition of the original system into the main system and a tagged server station as above. The further arguments for stochastic comparison of Markov processes are also easily extended for the case of an arbitrary $m\geq2$.
\end{proof}

\smallskip
From the above results for the Markov processes the statement of
the lemma follows. The stochastic inequalities between
$T_{m,n}(m-1)$ and $T_{m,0}(m-1)$ follow by the coupling
arguments. The lemma is completely proved.
\end{proof}

\section{Theorems on losses in $M/GI/m/n$ queueing
systems}\label{Main theorems}

The results obtained in the previous section enable us to
establish theorems for the number of losses in $M/GI/m/n$ queueing
systems during their busy periods.

\begin{thm}\label{Theorem4.1}
Under the assumption $\lambda=m\mu$, the expected number of losses
during a busy period of the $M/GI/m/n$ queueing system is the same
for all $n\geq1$.
\end{thm}

\begin{proof}
Consider the system $M/GI/m/n$ under the assumption $\lambda=m\mu$, and similarly to the construction in the proof of Lemma \ref{Lemma3.1} let us delete all the intervals where the number of customers in the system is less than $m$, and merge the corresponding ends.
 The process obtained is denoted
$\widehat{\mathbf{{y}}}_{m,n}(t)$. This is the 1-stationary
process of orbital busy periods.

The stationary departure process, together with the arrival
1-Poisson process of rate $\lambda$ and the number of waiting places
$n$ describes the stationary $M/G/1/n$ queue-length process (with
generally dependent service times). As soon as a busy period is
finished (in our case it is an orbital busy period, see Section
\ref{Markovian} for the definition), the system immediately starts a
new busy period by attaching a new customer into the system. This
unusual situation arises because of the construction of the process.
There are no idle periods, and servers all are continuously busy.
Thus, the busy period, which is considered here, is one of the busy
periods attached one after another.

Let $T$ be a large period of time, and during that time there are
$K(T)$ busy periods of the $M/G/1/n$ queueing system (which does not
contain idle times as mentioned). Let $L(T)$ and $\nu(T)$ denote the
number of lost and served customers during time $T$. We have the
formula:
\begin{equation}\label{4.6}
\lim_{T\to\infty}\frac{1}{{\rm E}K(T)}\Big({\rm E}L(T)+{\rm
E}\nu(T)\Big)=\lim_{T\to\infty}\frac{1}{{\rm E}K(T)}\Big(\lambda
T+{\rm E}K(T)\Big),
\end{equation}
the proof of which is given below.

Relationship \eqref{4.6} has the following explanation. The
left-hand side term ${\rm E}L(T)+{\rm E}\nu(T)$ is the expectation
of the number of lost customers plus the expectation of the number
of served customers during time $T$, and the right-hand side term
$\lambda T+{\rm E}K(T)$ is the expectation of the number of arrivals
during time $T$ plus the expected number of attached customers.

Relationship \eqref{4.6} can be proved by renewal arguments as
follows.

There are $m$ independent copies ${\bf x}_1^{(1)}(t)$, ${\bf
x}_1^{(2)}(t)$, \ldots, ${\bf x}_1^{(m)}(t)$ of the stationary
renewal process, which model the process $\widehat{\bf
y}_{m,n}(t)$. (In fact, we have $m$ 1-stationary processes, which
have the same distributions as $m$ independent renewal processes
${\bf x}_1^{(1)}(t)$, ${\bf x}_1^{(2)}(t)$, \ldots, ${\bf
x}_1^{(m)}(t)$.) Let $1\leq i\leq m$, and let $C_1$, $C_2$,\ldots
$C_{K_i(T)}$ be such points of busy period starts associated with
the renewal process $\mathbf{x}_1^{(i)}(t)$ (one of those $m$
independent and identically distributed renewal processes), where
$K_i(T)$ denotes the total number of these regeneration point
indexed by $i$. Denote also by $z_1$, $z_2$, \ldots, $z_{K_i(T)}$
the corresponding lengths of busy periods, by $\ell_1$,
$\ell_2$,\ldots,$\ell_{K_i(T)}$ the corresponding number of losses
during these $K_i$ busy periods, and by $n_1$,
$n_2$,\ldots,$n_{K_i(T)}$ the corresponding number of served
customers during these busy periods. Let
$T_i=z_1+z_2+\ldots+z_{K_i(T)}$, let
$L_i(T)=\ell_1+\ell_2+\ldots+\ell_{K_i(T)}$ and let
$\nu_i(T)=n_1+n_2+\ldots+n_{K_i(T)}$.

Since at the moments $C_1$, $C_2$,\ldots $C_{K_i(T)}$ of the busy
period starts the distribution of the above stationary Markov
process of residual times is the same, then the numbers of losses
$\ell_1$, $\ell_2$,\ldots,$\ell_{K_i(T)}$ and, respectively, the
numbers of served customers $n_1$, $n_2$,\ldots,$n_{K_i(T)}$
during each of these busy periods have the same distributions, and
one can apply the renewal reward theorem.

By the renewal reward theorem we have:
\begin{equation}\label{4.6+}
\lim_{T\to\infty}\frac{1}{m{\rm E}K_i(T)}\Big({\rm E}L_i(T)+{\rm
E}\nu_i(T)\Big)=\lim_{T\to\infty}\frac{1}{m{\rm
E}K_i(T)}\Big(\lambda \mathrm{E}T_i+m{\rm E}K_i(T)\Big).
\end{equation}
 Taking into account that
$$
\lim_{T\to\infty}\frac{\mathrm{E}K(T)}{\mathrm{E}K_i(T)}=m,
$$
$$
\lim_{T\to\infty}\frac{\mathrm{E}L(T)}{\mathrm{E}L_i(T)}=1,
$$
$$
\lim_{T\to\infty}\frac{\mathrm{E}\nu(T)}{\mathrm{E}\nu_i(T)}=1,
$$
and
$$
\lim_{T\to\infty}\frac{\mathrm{E}T_i}{T}=1,
$$
 because of the correspondence between the left- and right-hand
sides, from \eqref{4.6+} we arrive at \eqref{4.6}. Thus, we bypass
the fact that the times between departures are dependent, and thus
\eqref{4.6} is actually obtained by application of the renewal
reward theorem by a usual manner, like in the case of independent
times between departures (see e.g. Ross
\cite{Ross 2000}, Karlin and Taylor \cite{Karlin and Taylor 1975}).

Together with \eqref{4.6} we have
\begin{equation}\label{4.7}
\lim_{T\to\infty}\frac{1}{{\rm E}K(T)}{\rm
E}\nu(T)=\lim_{T\to\infty}\frac{1}{{\rm E}K(T)}m\mu T.
\end{equation}

Let us now introduce the following notation. Let $\zeta_n$ denote
the length of an orbital busy period, and let $L_n$ and $\nu_n$
correspondingly denote the number of lost and served customers
during that orbital busy period. Using the arguments of
\cite{Abramov 2001a}, we prove that ${\rm E}L_n=1$ for all
$n\geq1$.

Indeed, from (\ref{4.6}) and (\ref{4.7}) we have the equations:
\begin{equation}\label{4.8}
{\rm E}L_n+{\rm E}\nu_n=\lambda{\rm E}\zeta_n+1,
\end{equation}
\begin{equation}\label{4.9}
{\rm E}\nu_n=m\mu{\rm E}\zeta_n.
\end{equation}
The substitution $\lambda=m\mu$ into the system of equations
(\ref{4.8}) and (\ref{4.9}) yields ${\rm E}L_n=1$.

Hence, during an orbital busy period there is exactly one lost
customer in average for any $n\geq0$. To finish the proof we need in
a deeper analysis. First, we should find the expected number of
queueing periods during one orbital busy period. For this purpose,
one can use the similar construction by deleting all the intervals
when the number of customers in the system is not greater than $m$,
and merge the corresponding ends. The obtained process is denoted
$\widehat{\mathbf{y}}_{m+1,n}(t)$, and this is one stationary
process of queueing periods following one after another.

The structure of the process $\widehat{\mathbf{y}}_{m+1,n}(t)$ is
similar to that of the process $\widehat{\mathbf{y}}_{m,n}(t)$. The
process $\widehat{\mathbf{y}}_{m+1,n}(t)$ describes the stationary
$M/G/1/n-1$ queueing system, the service times of which are
generally dependent. As soon as one busy period in this system is
finished, a new customer starting a new busy period is immediately
attached into the system. Thus, the only difference between the
processes $\widehat{\mathbf{y}}_{m,n}(t)$ and
$\widehat{\mathbf{y}}_{m+1,n}(t)$ is that the numbers of waiting
places differ by value of parameter $n$. Therefore, using the
similar notation and arguments, one arrive at the conclusion that
the expected number of losses per queueing period is equal to 1 as
well. Therefore, in long-run period of time, the number of queueing
periods and orbital busy periods is the same in average. So, there
is exactly one queueing period per orbital busy period in average.

Therefore, during
a long-run period the number of events that the different Markov
processes $\widehat{\bf y}_{m-1,n}(t)$ change their values after
deleting queueing periods and merging the ends (as exactly
explained in the proof of Lemma \ref{Lemma3.1}) is the same in
average for all $n\geq1$,
and
the stationary characteristics of all of these Markov processes
$\widehat{\bf y}_{m-1,n}(t)$, given for different values $n$=1,2,\ldots, are the same. Hence,
the expectation $\mathrm{E}T_{m,n}(m-1)$ is the same for all
$n$=1,2,\ldots as well. (Recall that $T_{m,n}(m-1)$ denote the total time during a busy period when
$m-1$ servers are occupied.)

Hence, using Wald's identity connecting $\mathrm{E}T_{m,n}(m-1)$ with $\mathrm{E}L_{m,n}$ (the expected number of losses during a busy period) we arrive at
the desired result, since
$\mathrm{E}T_{m,n}(m-1)$ and the expectation of the number of
orbital busy periods during a busy period of $M/GI/m/n$ both are the same for all $n\geq1$.
\end{proof}

Application of Lemma \ref{Lemma3.1} and the arguments
of Theorem \ref{Theorem4.1} enables us to prove the following
result.

\begin{thm}
Let $\lambda=m\mu$. Then, under the assumption that $G(x)$ belongs
to the class NBU, for the number of losses in $M/GI/m/n$ queueing
systems, $n\geq1$, we have
\begin{equation}
\label{4.10}\mathrm{E}L_{m,n}=\frac{cm^m}{m!},
\end{equation}
where the constant $c\geq1$ depends on $m$ and the probability
distribution $G(x)$ but is independent of $n$.

Under the assumption that $G(x)$ belongs to the class NWU we have
\eqref{4.10} but with constant $c\leq1$.
\end{thm}

\begin{proof}
Notice first, that for the expected number of losses in $M/GI/m/0$
queueing systems we have
\begin{equation*}
\mathrm{E}L_{m,0}=\frac{m^m}{m!}
\end{equation*}
This result follows immediately from the Erlang-Sevastyanov formula
\cite{Sevastyanov 1957}, so that the expected number of losses
during a busy period of the $M/GI/m/0$ queueing system is the same
that this of the $M/M/m/0$ queueing system. The expected number of
losses during a busy period of the $M/M/m/0$ queueing system,
$\mathrm{E}L_{m,0}=\frac{m^m}{m!}$, is also derived in Section
\ref{Markovian}.

In the case where $G(x)$ belongs to the class NBU according to
Lemma \ref{Lemma3.1} we have
$\mathrm{E}T_{m,n}(m-1)\geq\mathrm{E}T_{m,0}(m-1)$, and therefore,
the expected number of orbital busy periods in the $M/GI/m/n$
queueing system ($n\geq1$) is not smaller than this in the
$M/GI/m/0$ queueing system. Therefore, repeating the proof of
Theorem \ref{Theorem4.1} leads to the inequality
$\mathrm{E}L_{m,n}\geq\mathrm{E}L_{m,0}$ and consequently to the
desired result. If $G(x)$ belongs to the class NWU, then we have
the opposite inequalities, and finally the corresponding result
stated in the formulation of the theorem.
\end{proof}

\section{Batch arrivals}
The case of batch arrivals is completely analogous to the case of
ordinary (non-batch) arrivals. In the case of a Markovian
$M^X/M/m/n$ queueing system one can also apply the level-crossing
method to obtain equations analogous to (\ref{2.6})-(\ref{2.14}).
The same arguments as in Sections \ref{Non-Markovian} and
\ref{Main theorems} in an extended form can be used for $M^X/GI/m/0$
queueing systems.



\end{document}